\documentclass{amsart}
\usepackage{amsmath}
\usepackage{amsthm}
\usepackage{amssymb}
\usepackage{amsfonts}
\usepackage{enumerate}
\usepackage{graphicx}

\normalsize

\newtheoremstyle{theoremstyle}
  {10pt}      
  {5pt}       
  {\itshape}  
  {}          
  {\bfseries} 
  {:}         
  {.5em}      
  {}          

\newtheoremstyle{examplestyle}
  {10pt}      
  {5pt}       
  {}          
  {}          
  {\bfseries} 
  {:}         
  {.5em}      
  {}          

\theoremstyle{theoremstyle}
\newtheorem{theorem}{Theorem}[section]
\newtheorem*{theorem*}{Theorem}

\newtheorem*{proposition*}{Proposition}

\newtheorem*{corollary*}{Corollary}

\theoremstyle{examplestyle}

\newtheorem{definition}[theorem]{Definition}
\newtheorem{definition*}{Definition}

\newtheorem{remark*}{Remark}

\newcommand{\comment}[1]{}

\newcommand{\Z}{\mathbb{Z}}
\newcommand{\R}{\mathbb{R}}

\newcommand{\pp}{0}

\begin{document}

\title{A Counterexample to King's Conjecture}


\author{Lutz Hille}

\address{Fakult\"at f\"ur Mathematik, Universit\"at Bielefeld, Postfach 100131,
33501 Bielefeld, Germany}

\author{Markus Perling}

\address{Institut Fourier - UMR5582, 100 rue des Maths, BP 74, 38402, St. Martin
d'Heres, France}

\thanks{The second author is supported by a grant of the German Research Council (DFG)}

\email{hille@math.uni-bielefeld.de}

\email{perling@mozart.ujf-grenoble.fr}

\date{January 2006}

\begin{abstract}
King's conjecture states that on every smooth complete toric variety $X$ there exists a
strongly exceptional collection which generates the bounded derived category
of $X$ and which consists of line bundles. We give a counterexample to this
conjecture. This example is just the Hirzebruch surface $\mathbb{F}_2$ iteratively
blown up three times, and we show by explicit computation
of cohomology vanishing that there exist no strongly exceptional sequences of
length 7.
\end{abstract}

\maketitle

\section{Introduction}

It is a widely open question whether on a given smooth algebraic variety $X$ (say,
complete and smooth), there exists a {\em tilting sheaf}. A tilting sheaf is a sheaf
$\mathcal{T}$ which generates the bounded derived category $\mathcal{D}^b(X)$ of $X$
and
$\operatorname{Ext}^k(\mathcal{T}, \mathcal{T}) = 0$ for all $k > 0$. For such
$\mathcal{T}$, the functor
\begin{equation*}
\operatorname{{\bf R}Hom}(\mathcal{T}, \, . \,) : D^b(X) \longrightarrow
D^b(A-\operatorname{mod}),
\end{equation*}
where $A := \operatorname{End}(\mathcal{T})$ is the endomorphism algebra, induces an
equivalence of categories (see \cite{Rudakov90}, \cite{Bondal90}, \cite{Beilinson78}).
The existence of a tilting sheaf implies that the Grothendieck group of $X$ is finitely
generated and free, so that in general such sheaves can not exist. However,
so far there are a number of positive examples known, including projective spaces,
del Pezzos, certain homogeneous spaces, and some higher dimensional Fanos.
An obvious testbed for the existence of tilting
sheaves are the toric varieties. There is a quite strong conjecture which was first
stated by King:\\

\nocite{AKO}

{\bf Conjecture} \cite{King2} {\bf :} Let $X$ be a smooth complete
toric variety. Then $X$ has a tilting sheaf which is a direct sum of line bundles.\\

If a tilting sheaf decomposes into a direct sum of line bundles, its direct summands
$\mathcal{T} = \bigoplus_{i = 1}^t \mathcal{L}_i$ form a so-called {\em strongly
exceptional sequence}, i.e. $\operatorname{Ext}^k(\mathcal{L}_i, \mathcal{L}_j) = 0$
for all $i, j$ and all $k > 0$, and  --- after eventually reordering the
$\mathcal{L}_i$ --- $\operatorname{Hom}(\mathcal{L}_i, \mathcal{L}_j) = 0$ for $i >
j$. Moreover, $t$ is the rank of the Grothendieck group of $X$.

\comment{If $\mathcal{T}$ splits into a direct sum of line bundles, i.e. $\mathcal{T} =
\bigoplus_{i = 1}^t \mathcal{L}_i$, this implies for the summands that we have
$\operatorname{Ext}^k(\mathcal{L}_i, \mathcal{L}_j) = 0$ for all $i, j$ and all
$k > 0$. Thus the sequence of line bundles $\mathcal{L}_1, \dots, \mathcal{L}_t$ forms
a so-called {\em strongly exceptional sequence}. Moreover, the condition that
$\mathcal{T}$ generates $\mathcal{D}^b(X)$ implies that $t$ equals the rank of the
Grothendieck group of $X$.}

It would be very nice if there existed easy-computable tilting sheaves on toric
varieties, and indeed there are known a lot of positive examples in favor of the
conjecture
(see \cite{CostaMiroRoig}, \cite{Kawamata1}, \cite{Hille1},
\cite{CrawSmith}). Computer experiments also look promising in many directions.
However, the conjecture remained somewhat mysterious so far and,
as it turns out, it is false in general. It is the purpose of this paper to present a
counterexample.

\nocite{BergmanProudfoot}

\begin{figure}[ht]
\includegraphics[height=6cm,width=6cm]{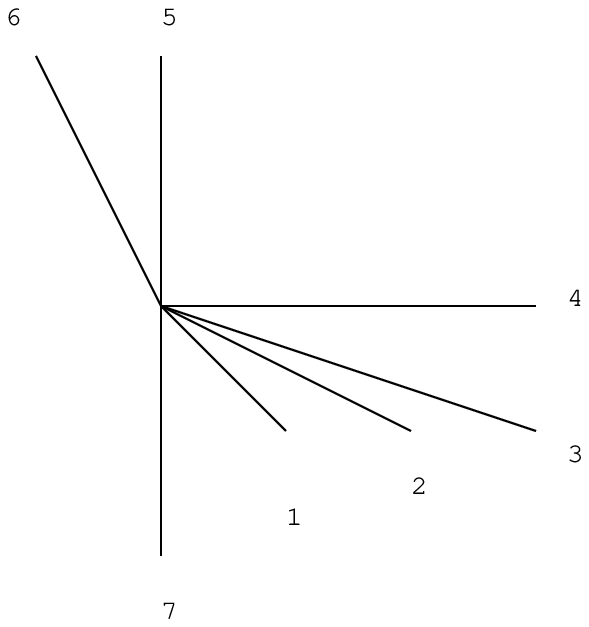}\quad
\caption{The fan}\label{fan}
\end{figure}

Our counterexample is the toric surface $X$ as shown in figure \ref{fan}, which can
be obtained by iteratively blowing up the Hirzebruch surface $\mathbb{F}_2$ three
times. In coordinates, the primitive vectors
of its rays are given by $l_1 = (1, -1)$, $l_2 = (2, -1)$, $l_3 = (3, -1)$, $l_4 =
(1, 0)$, $l_5 = (0, 1)$, $l_6 = (-1, 2)$, $l_7 = (0, -1)$.
Note that the rank of the Grothendieck group of $X$ is 7.

To show that there do not exist any strongly exceptional sequences of length 7 on this
surface, we will perform explicit computations in the Picard group to determine
cohomology
vanishing. More precisely, note that if $\mathcal{L}_1, \dots, \mathcal{L}_t$ is a
strongly exceptional sequence, then also $\mathcal{L}_1 \otimes \mathcal{L}', \dots,
\mathcal{L}_t \otimes \mathcal{L}'$ is strongly exact, where $\mathcal{L}'$ is any
line bundle. So one can assume without loss of generality that the sequence contains
the structure sheaf. Then a necessary condition for the bundles in the sequence is that
all the higher cohomology groups of the bundles and of their dual bundles vanish,
i.e. $H^k(X, \mathcal{L}_i) = H^k(X, \mathcal{L}^*_i) = 0$ for all $i$
and all $k > 0$. This is a rather strong condition on the sheaves and our main
computation will be to compile a complete list of such bundles for our surface $X$.
After having obtained this classification, we deduce by simple inspection that a
strongly exceptional sequence of length 7 and consisting of line bundles does not exist.

{\bf Overview:} In section \ref{setup} we state everything we need to know about
cohomology of line bundles on toric surfaces and we describe in more detail our method
of computation. In section \ref{classification} all bundles are classified which have
the property that the higher cohomologies of the bundles themselves and of their dual
bundles vanish. In section \ref{listsection} we present the complete classification
obtained in section \ref{classification} and we show by inspection that there exist no
strongly exceptional sequences of length 7 on $X$.

\section{The setup}\label{setup}

In this section we recall basic facts on cohomology of line bundles on a toric
surface and we describe our method of computation. For general information about
toric varieties we refer to the books \cite{Oda}, \cite{Fulton}.

\subsection{Generalities on toric line bundles}
Let $X$ be a complete smooth toric surface on which the torus $T$ acts. The variety $X$
is described by a fan $\Delta$ which is contained in a $2$-dimensional vector space
$N_\R := N \otimes_\Z \R$, where $N \cong \Z^2$. We denote by $\Delta(1)$
the set of rays, that is, of one-dimensional cones of $\Delta$. As $X$ is a complete
surface, the fan is completely determined by the rays. We denote the rays by
$\rho_1, \dots, \rho_n$, enumerated in counterclockwise order, and $l_1, \dots,
l_n$ the primitive vectors of the rays. To any $\rho_i$ there is associated a
$T$-invariant divisor $D_i$, and every divisor $D$ can, up to rational equivalence,
written as a sum of these invariant divisors, i.e. $D = \sum_{i = 1}^n c_i D_i$.
We denote $M \cong \Z^2$ the character group of the torus acting on $X$ and we set
$M_\R := M \otimes_\Z \R$. The lattice $N$ is in a natural way dual to $M$, and
the primitive vectors $l_i$ are integral linear forms on $M$ (and on $M_\R$,
respectively). There is a short exact sequence
\begin{equation*}
0 \longrightarrow M \overset{A}{\longrightarrow} \Z^{\Delta(1)} \longrightarrow
\operatorname{Pic}(X) \longrightarrow 0,
\end{equation*}
where the matrix $A$ is composed of the $l_i$ as row vectors. This sequence is split
exact. More precisely, if we choose two of the $l_i$, for instance $l_{n - 1}$ and
$l_n$, which form a $\mathbb{Z}$-basis of $N$, then the divisors $D_1, \dots,
D_{n - 2}$ form a $\Z$-basis of $\operatorname{Pic}(X)$. So every divisor $D$ has
a unique representation $D = \sum_{i = 1}^{n - 2} c_i D_i$.

Now let $D = \sum_{i = 1}^n c_i D_i$ be any $T$-invariant divisor. $D$ in a natural
way defines an affine hyperplane arrangement $\mathcal{H}_D =  \{H_1, \dots, H_n\}$
in the vector space $M_\R$, where
\begin{equation*}
H_i = \{m \in M_\R \mid l_i(m) = - c_i\}.
\end{equation*}
All information on the cohomology of the line bundle $\mathcal{O}(D)$ is contained in
the chamber structure $\mathcal{H}_D$ (or more precisely, in the intersection of this
chamber structure with the lattice $M$). Recall that the $T$-action induces an
eigenspace decomposition on the cohomology groups of $\mathcal{O}(D)$:
\begin{equation*}
H^i\big(X, \mathcal{O}(D)\big) = \bigoplus_{m \in M} H^i\big(X,\mathcal{O}(D)\big)_m.
\end{equation*}
The dimension of $H^i\big(X,\mathcal{O}(D)\big)_m$ as a $k$-vector space is determined
by the {\em signature} of $m$ with respect to the arrangement $\mathcal{H}_D$:

\begin{definition}
Let $D = \sum_{i = 1}^n c_i D_i$ be a $T$-invariant divisor on $X$. Then for every
$i = 1, \dots, n$ we define a signature
\begin{equation*}
\Sigma^D_i : M \longrightarrow \{+, -, \pp\},
\end{equation*}
where $\Sigma^D_i(m) = +$ if $l_i(m) > - c_i$, $\Sigma^D_i(m) = -$ if $l_i(m) < - c_i$
and $\Sigma^D_i(m) = \pp$ if $l_i(m) = - c_i$. Moreover, we denote
\begin{equation*}
\Sigma^D : M \longrightarrow \{+, -, \pp\}^n,
\end{equation*}
where $\Sigma^D(m)$ is the tuple $\big(\Sigma^D_1(m), \dots, \Sigma^D_n(m)\big)$.
\end{definition}

Below we will mostly work with only one $D$ at a time, which will be clear from the
context. So usually we will omit the reference to $D$ in the notation, i.e. we will
mostly write $\Sigma(m)$ instead of $\Sigma^D(m)$.

Given the signature $\Sigma^D(m)$, the computation of $H^i\big(X, \mathcal{O}(D)
\big)_m$ is straightforward. For $H^2$, we have:
\begin{equation*}
\dim H^2\big(X, \mathcal{O}(D) \big)_m =
\begin{cases}
1 & \text{ if } \Sigma^D(m) = \{-\}^n \\
0 & \text{ else}.
\end{cases}
\end{equation*}
For $H^1$, we have to consider the {\em $-$-intervals}. For a given signature
$\Sigma^D(m)$, a $-$-interval is a connected sequence of $-$ with respect to the
circular order of the $\rho_i$. For example, assume that $\Delta(1)$ consists of
7 elements enumerated in circular order. Then the signature $+ - - + + - +$ has
two $-$-intervals. Note that due to the {\em circular} ordering of the rays, the
signature $- - + + + - -$ has only one $-$-interval. We have:
\begin{equation*}
\dim H^1\big(X, \mathcal{O}(D) \big)_m =
\begin{cases}
\text{ the number of $-$-intervals } - 1 \text{ if there exists at least one $-$-interval}, \\
0 \text{ else}.
\end{cases}
\end{equation*}
Thus $H^1\big(X, \mathcal{O}(D)\big)$ vanishes if and only if the signatures
$\Sigma^D(m)$, as $m$ runs through $M$, have at most one $-$-interval.

\subsection{Method of computation}

Let $\mathcal{L}_1, \dots, \mathcal{L}_t$ be a strongly exceptional sequence of line
bundles, i.e. we have $\operatorname{Ext}^k(\mathcal{L}_i, \mathcal{L}_j) = 0$ for
all $i, j$ and all $k > 0$. There is a natural isomorphism
$\operatorname{Ext}^k(\mathcal{L}_i, \mathcal{L}_j) \cong H^k(X, \mathcal{L}_i^*
\otimes \mathcal{L}_j)$, where $\mathcal{L}_i^* = \operatorname{Hom}(\mathcal{L}_i,
\mathcal{O}_X)$ denotes the dual bundle.
 By this we can assume without loss of generality that one of
the $\mathcal{L}_i$ is just the structure sheaf $\mathcal{O}_X$, i.e.
$\mathcal{L}_1, \dots, \mathcal{L}_t$ is a strongly exceptional sequence if and only
if $\mathcal{L}_i^* \otimes \mathcal{L}_1, \dots, \mathcal{L}_i^* \otimes
\mathcal{L}_t$ is a strongly exceptional sequence. If $\mathcal{O}_X$ is part of
the sequence, this in turn implies a rather strong condition on the cohomologies of
the other bundles. Namely, for every $\mathcal{L}_i$ we have:
\begin{equation*}
H^k(X, \mathcal{L}_i) = H^k(X, \mathcal{L}_i^*) = 0 \text{ for all } k > 0.
\end{equation*}
Thus, to show that our toric surface does not have a strongly exceptional sequence of
length 7, we proceed in 2 steps:
\begin{enumerate}[(i)]
\item We classify all line bundles where higher cohomologies of the bundle itself as
well as of its dual vanish. It turns out that the list of such bundles has a rather
short description, although it is not finite.
\item After having obtained the list, we show by exclusion that there are no strongly
exceptional sequences of length 7.
\end{enumerate}

\begin{figure}[ht]
\

\includegraphics[height=10cm,width=10cm]{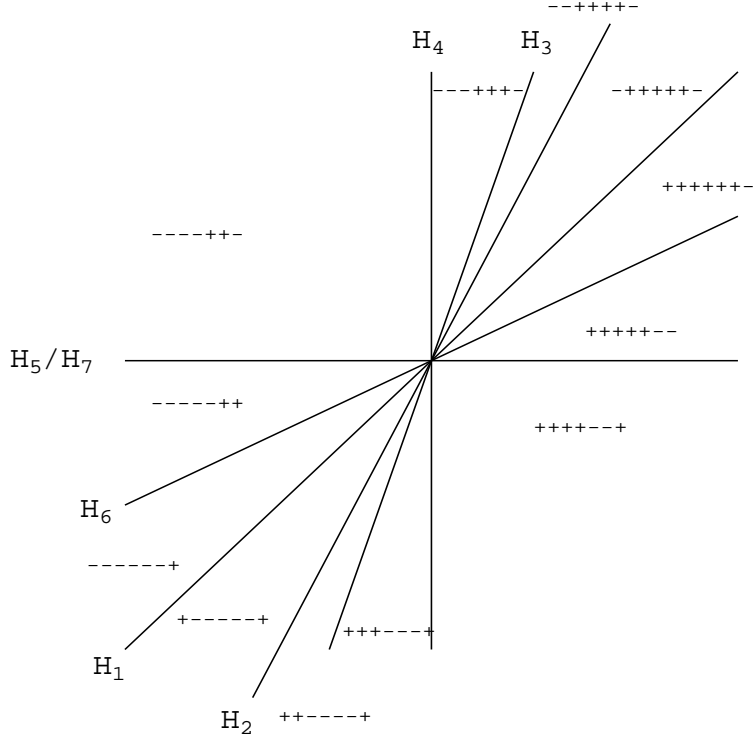}\qquad
\caption{The central arrangement}\label{centralarrangement}
\end{figure}

\begin{figure}[ht]
\
\includegraphics[height=10cm,width=10cm]{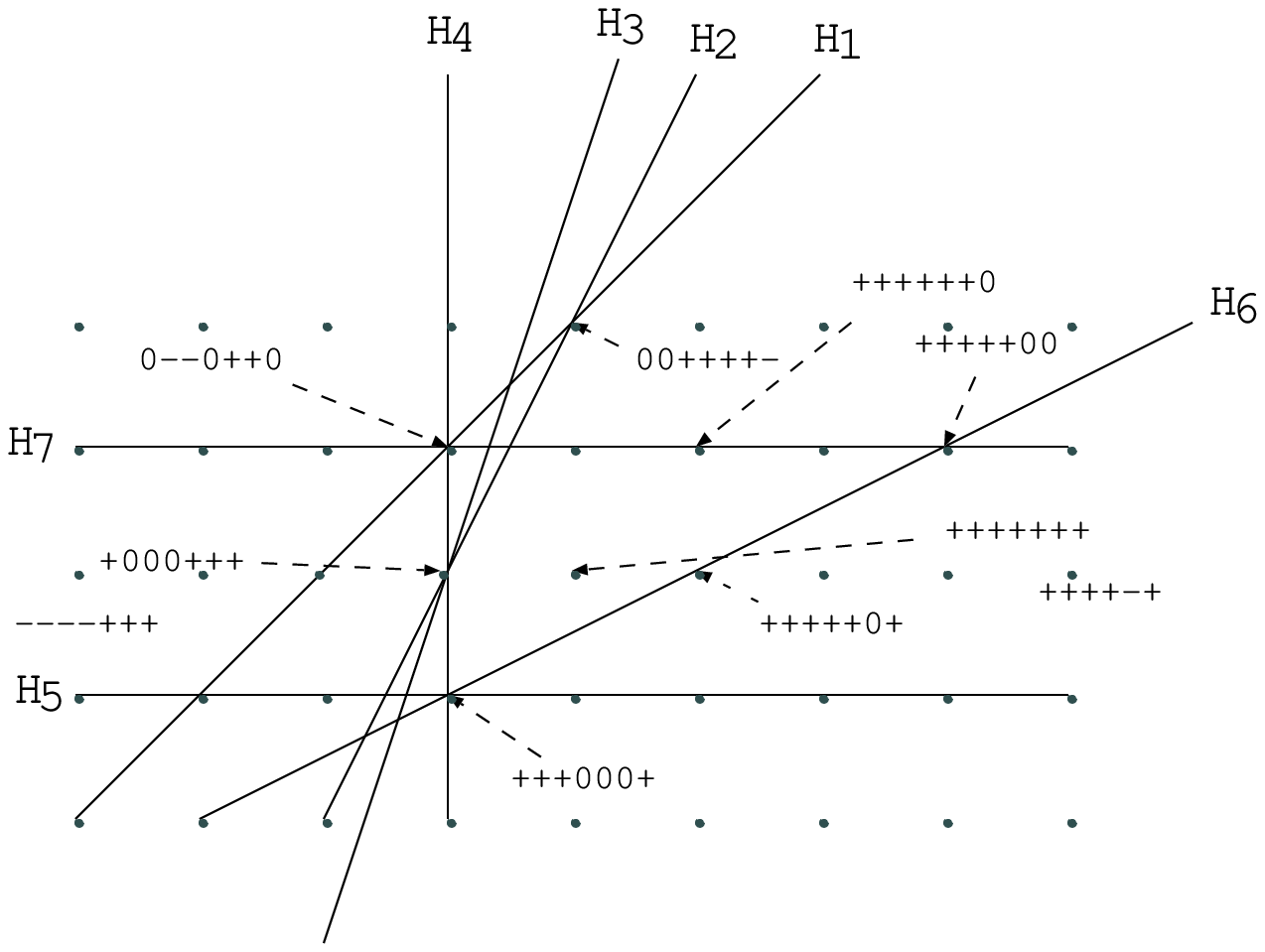}\quad
\caption{A deformation of the central arrangement}\label{deformedarrangement}
\end{figure}

Figure \ref{centralarrangement} shows the arrangement which belongs
to the structure sheaf. We see that this arrangement is central and induces a chamber
decomposition of the space $M_\R$, consisting of unbounded chambers. To every
chamber there is associated a signature which we have indicated in the picture. Note
that in fact there are some more signatures which are not shown. For instance,
the points lying on the line between the chambers with signatures $++++++-$ and
$-+++++-$ have signature $\pp+++++-$. The origin has signature $\pp\pp\pp\pp\pp
\pp\pp$. Figure \ref{deformedarrangement} shows a deformation of this
central arrangement which belongs to the divisor $D = - (4 D_1 + 7 D_2 + 11 D_3 +
4 D_4 + 2 D_5)$.


 As we can see, moving the hyperplanes creates new chambers with new
signatures. There are two new unbounded chambers with signatures $----+++$ and
$+++++-+$, respectively, which obviously have no influence on the comohology of
$\mathcal{O}(D)$. The other chambers are all bounded and thus contain only a finite
number of lattice points (i.e. points in $M$). We have indicated the signatures of
some of these points in the picture. As one can check, most of these signatures give
not rise to nonvanishing cohomology, the only exception being the point with
signature $+++++++$. Recall that we are interested in the classification of line
bundles which have no higher cohomology and whose duals have also no higher cohomology.
So, if
there is an inequality $l_i(m) < - c_i$ (or $l_i(m)> - c_i$, respectively), then we
have
$l_i(-m) > c_i$ ($l_i(-m) < c_i$, respectively), whereas for $l_i(m) = -c_i$ we have
$l_i(-m) = c_i$. In our example the
signature of $m$ with $\Sigma^D(m) = +++++++$ becomes $\Sigma^{-D}(-m) = -------$ for
the dual bundle, which therefore has nonvanishing $H^2$.

\begin{figure}[ht]
\includegraphics[height=5cm,width=9cm]{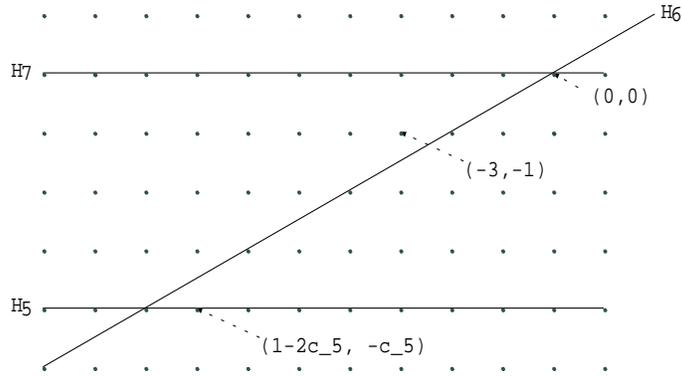}\quad
\caption{A partial arrangement}\label{partialarrangement}
\end{figure}

We give one more example and some more notation. In many situation it will not be
necessary to know the complete signature of some point $m \in M$. Therefore we
define:
\begin{definition}
A {\em partial} signature is given by
\begin{equation*}
\Sigma^D : M \longrightarrow \{+, -, \pp, *\}^n
\end{equation*}
which is a signature for some subset $I$ of $\{1, \dots, n\}$ such that
$\big(\Sigma^D(m)\big)_i = \Sigma^D_i(m)$ for
$i \in I$ and $\big(\Sigma^D(m)\big)_i = *$ for $i \notin I$.
\end{definition}
For us it is convenient to use the same symbol for signatures and partial signatures.
Let us give an explicit example. Assume $D = \sum_{i = 1}^5 c_i D_i$ and
$c_5 > 0$. Now consider the point $m$ in $M$ which has the coordinates $(1 - 2 c_5,
-c_5)$ (see figure \ref{partialarrangement}). Its partial signature with respect to
the linear forms $l_5, l_6, l_7$ is $\Sigma(m) = ****\pp-+$. Our aim is to derive
conditions on the values of the $c_i$. Evidently, any complete
signature which is obtained by filling the $*$'s has at least one $-$-interval.
Moreover, if any of the $*$'s becomes a $-$, the signature has at least two
$-$-intervals, and any corresponding line bundle will have nonvanishing $H^1$. So,
a necessary condition is that $\Sigma^D_i(m) \in \{+, \pp\}$ for $i = 1, \dots, 4$
and any valid divisor $D$. This in turn implies:
\begin{equation}\label{firstgeneralconditions}
\begin{split}
c_1 & \geq c_5 - 1 \\
c_2 & \geq 3 c_5 - 2\\
c_3 & \geq 5 c_5 - 3\\
c_4 & \geq 2 c_5 - 1.
\end{split}
\end{equation}
Now the point $(-3, -1)$ has partial signature $\Sigma(-3,-1) = ****+++$, and the
above conditions on $c_1, \dots, c_4$ imply that for $c_5 > 3$ this point always
has signature $+++++++$, and thus we have nonvanishing $H^2$. Hence, we conclude
$c_5 \leq 3$.

\section{Classification of line bundles without higher cohomology}\label{classification}

In this section we do the complete classification of line bundles for our toric
surface which have the property that the higher cohomologies vanish for both, the
bundle itself and its dual. As explained in the previous section, we can always
assume that a line bundle $\mathcal{L}$ is uniquely represented by an invariant
divisor $D = \sum_{i = 1}^5 c_i D_i$, and every tuple of numbers $(c_1, \dots,
c_5)$ represents a unique isomorphism class in $\operatorname{Pic}(X)$. As we already
have seen, a necessary condition is that $c_5 \leq 3$. Moreover, as it does not matter
if we deal with a bundle or its dual, we can assume without loss of generality that
$c_5 \geq 0$. So, this leaves us with four possible values for $c_5$. Our
classification will be done by subsequential case distinctions which on the toplevel
are guided by the four possible values of $c_5$.

Note that in the sequel for a given bundle we will use phrases like
``has cohomology'' if either the bundle itself or its dual
has a nonvanishing higher cohomology group.

\subsection{$c_5 = 3$}

Recall that the partial signature of the point $(-3, -1)$ is $\Sigma^D(-3, -1) =
****+++$. By the conditions (\ref{firstgeneralconditions}), we immediately obtain
the partial signature $*++++++$. So, the only way to prevent $H^2$ to show up in the
dual bundle is $\Sigma_1^D(-3,-1) = \pp$ (then, for the dual bundle, we have the
signature $\Sigma^{-D}(3, 1) = \pp------$). This in turn means that $c_1 = 2$.
But then, we have $\Sigma^D(-4, -2) =
\pp++++\pp+$, and thus $\Sigma^{-D}(4, 2) = \pp----\pp-$, hence we get $H^1$. We
conclude that there are no divisors with $c_5 = 3$ and vanishing cohomology.

\subsection{$c_5 = 2$}
Here, conditions (\ref{firstgeneralconditions}) read:
\begin{align*}
c_1 & \geq 1, \\
c_2 & \geq 4, \\
c_3 & \geq 7, \\
c_4 & \geq 3.
\end{align*}
We first consider $c_4 = 3$. Then $\Sigma(-3, 0) = ***\pp++\pp$. If one of the $*$'s
is replaced by $+$, this implies that the dual signature will have at least two
$-$-intervals, independent on the other substitutions. So we obtain:
\begin{align*}
c_1 & \leq 3, \\
c_2 & \leq 6, \\
c_3 & \leq 9.
\end{align*}
We treat these 27 possibilities case by case. First, let $c_1 = 1$. Then we have
$\Sigma(-2, -1) = \pp**++\pp+$, and so we have $H^1$, leaving only 18 more cases.
For these we write a table:
\begin{center}
\begin{tabular}{|c|c|c|c|c|} \hline
$c_1$ & $c_2$ & $c_3$ & $m$ & $\Sigma(m)$\\ \hline
2 & 4 & 8 & $(-3, -2)$ & $+\pp+\pp\pp-+$\\ \hline
2 & 6 & --- & $(-3, -1)$ & $\pp+*\pp+++$\\ \hline
2 & --- & 9 & $(-3, -1)$ & $\pp*+\pp+++$\\ \hline
3 & 4 & --- & $(-2, 0)$ & $+\pp*+++\pp$\\ \hline
3 & 5 & 9 & $(-3, -1)$ & $+\pp+\pp+++$\\ \hline
3 & 6 & 7 & $(-2, 1)$ & $\pp+\pp+++-$\\ \hline
\end{tabular}
\end{center}

This table contains a list of all values which have cohomology. For given values of
$c_1, c_2, c_3$, the fourth columns contains a lattice point $m \in M$ with bad
signature, which is displayed in the fifth column. Sometimes it suffices to display
only a partial signature. Then the box for the corresponding $c_i$ contains a dash
(---). All tuples which are not displayed in the above table represent cohomology
free line bundles, namely:

\begin{center}
\begin{tabular}{|c|c|c|c|c|} \hline
$c_1$ & $c_2$ & $c_3$ & $c_4$ & $c_5$\\ \hline
2 & 4 & 7 & 3 & 2 \\ \hline
2 & 5 & 7 & 3 & 2 \\ \hline
2 & 5 & 8 & 3 & 2 \\ \hline
3 & 5 & 7 & 3 & 2 \\ \hline
3 & 5 & 8 & 3 & 2 \\ \hline
3 & 6 & 8 & 3 & 2 \\ \hline
3 & 6 & 9 & 3 & 2 \\ \hline
\end{tabular}
\end{center}

Now for $c_4 = 4$. We have $\Sigma(-3, -2) = ***+\pp-+$ and thus we get the bounds
\begin{align*}
c_1 & \geq 2, \\
c_2 & \geq 5, \\
c_3 & \geq 8.
\end{align*}
Moreover, we have $\Sigma(-4, 0) = ***\pp++\pp$ and thus
\begin{align*}
c_1 & \leq 4, \\
c_2 & \leq 8, \\
c_3 & \leq 12.
\end{align*}
Further, $\Sigma(-3, -1) = ***++++$ and hence
\begin{align*}
c_1 & \leq 2 \text{ or} \\
c_2 & \leq 5 \text{ or} \\
c_3 & \leq 8.
\end{align*}
We have $\Sigma(-4, -3) = ***\pp--+$ which implies
\begin{equation*}
c_3 \geq 9
\end{equation*}
and so the case $c_3 \leq 8$ cannot occur. Also, the conditions imply that
either $c_1 = 2$ or $c_2 = 5$, thus leaving 24 possibilities.

We first consider the case $c_1 = 2$. Then we have $\Sigma(-3, -2) = \pp**\pp\pp\pp+$, which implies $c_2 \leq 6$ and $c_3 \leq 10$. Now we take $c_2 = 5$. Then
$\Sigma(-4,-3) = +\pp*\pp--+$, so that we must have $c_3 = 9$.
For $c_2 = 6$, we have $\Sigma(-4, -2) = \pp\pp*\pp\pp\pp+$, which implies $c_3 = 10$.
Indeed, we have found:

\begin{center}
\begin{tabular}{|c|c|c|c|c|} \hline
$c_1$ & $c_2$ & $c_3$ & $c_4$ & $c_5$\\ \hline
2 & 5 & 9 & 4 & 2 \\ \hline
2 & 6 & 10 & 4 & 2 \\ \hline
\end{tabular}
\end{center}

Now we consider $c_2 = 5$. We can assume that $c_1 \geq 3$. Assume that $c_3 \geq 10$.
Then $\Sigma(-4, -3) = +\pp+\pp--+$, so we have cohomology, hence $c_3 = 9$. For
$c_1 = 4$, we have $\Sigma(-4, -3) = + \pp+\pp--+$ and thus cohomology, hence
$c_1 = 3$, and indeed we have found:

\begin{center}
\begin{tabular}{|c|c|c|c|c|} \hline
$c_1$ & $c_2$ & $c_3$ & $c_4$ & $c_5$\\ \hline
3 & 5 & 9 & 4 & 2 \\ \hline
\end{tabular}
\end{center}

Now we go on with $c_4 \geq 5$. Then we have $\Sigma(-4, -2) = ***+\pp\pp+$, which
yields the conditions
\begin{align*}
c_1 & \geq 4, \\
c_2 & \geq 7, \\
c_3 & \geq 11.
\end{align*}
The signature $\Sigma(-3, -1)$ as before implies
\begin{align*}
c_1 & \leq 2 \text{ or} \\
c_2 & \leq 5 \text{ or} \\
c_3 & \leq 8.
\end{align*}
both conditions cannot be fulfilled simultaneously, and hence, for $c_4 \geq 5$ there
are no cohomologyfree bundles.

\subsection{$c_5 = 1$}
Again, we start with the conditions (\ref{firstgeneralconditions}), which read
\begin{align*}
c_1 & \geq 0, \\
c_2 & \geq 1, \\
c_3 & \geq 2, \\
c_4 & \geq 1.
\end{align*}
Now we go for the different cases for $c_1$.

$c_1 = 0$. Then we have $\Sigma(-1, -1) = \pp***\pp-+$ so that all of the $*$'s can
only be substituted by $\pp$'s, and thus $c_2 = 1$, $c_3 = 2$, $c_4 = 1$ and indeed we
have found:

\begin{center}
\begin{tabular}{|c|c|c|c|c|} \hline
$c_1$ & $c_2$ & $c_3$ & $c_4$ & $c_5$\\ \hline
0 & 1 & 2 & 1 & 1 \\ \hline
\end{tabular}
\end{center}
with no other possibilities left.

$c_1 = 1$. We have $\Sigma(-2, -1) = \pp***\pp\pp+$ which implies
\begin{align*}
c_2 & \leq 3, \\
c_3 & \leq 5, \\
c_4 & \leq 2.
\end{align*}
Let $c_4 = 1$, then $\Sigma(-1, 0) = \pp**\pp++\pp$, which implies
\begin{align*}
c_2 & \leq 2, \\
c_3 & \leq 3.
\end{align*}
From these four cases, only $c_2 = 1, c_3 = 2$ has cohomology, as in this case
$\Sigma(-1, -1) = +\pp+\pp\pp-+$. We have found:
\begin{center}
\begin{tabular}{|c|c|c|c|c|} \hline
$c_1$ & $c_2$ & $c_3$ & $c_4$ & $c_5$\\ \hline
1 & 1 & 2 & 1 & 1 \\ \hline
1 & 2 & 2 & 1 & 1 \\ \hline
1 & 2 & 3 & 1 & 1 \\ \hline
\end{tabular}
\end{center}

Now let $c_4 = 2$. Then $\Sigma(-2, -1) = +**+\pp-\pp$, so that
\begin{align*}
c_2 & \geq 2, \\
c_3 & \geq 3,
\end{align*}
leaving six cases.
We write a table as before:
\begin{center}
\begin{tabular}{|c|c|c|c|} \hline
$c_2$ & $c_3$ & $m$ & $\Sigma(m)$\\ \hline
2 & 3 & $(-3, -2)$ & $+\pp-\pp--+$\\ \hline
2 & 5 & $(-3, -2)$ & $+\pp+\pp--+$\\ \hline
3 & 3 & $(-1, 0)$ & $\pp+\pp+++\pp$\\ \hline
\end{tabular}
\end{center}
and thus we have found:
\begin{center}
\begin{tabular}{|c|c|c|c|c|} \hline
$c_1$ & $c_2$ & $c_3$ & $c_4$ & $c_5$\\ \hline
1 & 2 & 4 & 2 & 1 \\ \hline
1 & 3 & 4 & 2 & 1 \\ \hline
1 & 3 & 5 & 2 & 1 \\ \hline
\end{tabular}
\end{center}

$c_1 \geq 2$.
Now for any $c_1 \geq 2$, the point $(1 - c_1, 0)$ has signature $\Sigma(1 - c_1, 0) =
+***++\pp$. So, we obtain general conditions
\begin{equation*}
\begin{split}
c_2 & \geq 2 c_1 - 1, \\
c_3 & \geq 3 c_1 - 2, \\
c_4 & \geq c_1.
\end{split}
\end{equation*}
We obtain another general condition as follows. Consider the signature
$\Sigma(-c_1 - 1, 0) = -***++\pp$. Assume that $c_4 \geq c_1 + 2$. Then
$\Sigma(-c_1, 0) = -+**++\pp$ and so the $*$'s can only be replaced by $+$'s, hence
$c_2 \geq 2 c_1 + 3$. The signature $\Sigma(-c_2 - 2, - 1)$ then becomes either
$-\pp*-\pp+\pp$ or $-+*-\pp+\pp$ which both are bad. Thus $c_4$ must be strictly
smaller than $c_1 + 2$, and we have:
\begin{equation*}
c_4 \in \{c_1, c_1 + 1\} \text{ for any value of $c_1 \geq 2$}.
\end{equation*}
Now consider the signature $\Sigma(-c_1 - 1, -1) = \pp***\pp++$, which yields the
following restrictions:
\begin{equation*}\label{fourthgeneralconditions}
\begin{split}
c_2 & \leq 2 c_1 + 1, \\
c_3 & \leq 3 c_1 + 2. \\
\end{split}
\end{equation*}
Now assume that $c_4 = c_1$. From the signature $\Sigma(-c_1, 0) = \pp**\pp++\pp$ we
get immediately the conditions
\begin{equation*}
\begin{split}
c_2 & \leq 2 c_1, \\
c_3 & \leq 3 c_1. \\
\end{split}
\end{equation*}
If $c_2 = 2 c_1 - 1$, we have the signature $\Sigma(-c_1, -1) = +\pp*\pp\pp\pp +$,
respectively $\Sigma(-2, -1) = +\pp*\pp\pp++$, for the case $c_1 = 2$. In either case,
we get:
\begin{equation*}
c_3 \leq 3 c_1 - 1.
\end{equation*}
For $c_2 = 2 c_1$, we have the signature $\Sigma(1 - c_1, 1) = \pp+*+++-$, hence the
$*$ cannot be replaced by $-$ or $\pp$, thus we get $c_3 \geq 2 c_1 - 1$. We cannot
find any more restrictions and in fact we have found inifinite series of
cohomology-free line bundles:
\begin{center}
\begin{tabular}{|c|c|c|c|c|} \hline
$c_1$ & $c_2$ & $c_3$ & $c_4$ & $c_5$\\ \hline
$k \geq 2$ & $2k - 1$ & $3k - 2$ & $k$ & 1 \\ \hline
$k \geq 2$ & $2k - 1$ & $3k - 1$ & $k$ & 1 \\ \hline
$k \geq 2$ & $2k$ & $3k - 1$ & $k$ & 1 \\ \hline
$k \geq 2$ & $2k$ & $3k$ & $k$ & 1 \\ \hline
\end{tabular}
\end{center}

Now let $c_4 = c_1 + 1$. The signature $\Sigma(-c_4, -1) = +**\pp--+$ yields
\begin{equation*}
\begin{split}
c_2 & \geq 2 c_1, \\
c_3 & \geq 3 c_1 + 1.
\end{split}
\end{equation*}
This leaves four possibilities of which we can only exclude the case $c_2 = 2 c_1$,
$c_2 = 3 c_1 + 2$. Here we distinguish cases $c_1 = 2, 3, \geq 4$. For $c_1 = 2$, we
have $\Sigma(-3, -2) = +\pp+\pp--+$, for $c_1 = 3$ we have $\Sigma(-4, -2) =
+\pp+\pp\pp-+$ and for $c_1 \geq 4$ we have $\Sigma(-c_1 - 1, -2) = +\pp+\pp+-+$, all
of which are bad signatures. So, we have extracted three more series:
\begin{center}
\begin{tabular}{|c|c|c|c|c|} \hline
$c_1$ & $c_2$ & $c_3$ & $c_4$ & $c_5$\\ \hline
$k \geq 2$ & $2k$ & $3k + 1$ & $k + 1$ & 1 \\ \hline
$k \geq 2$ & $2k + 1$ & $3k + 1$ & $k + 1$ & 1 \\ \hline
$k \geq 2$ & $2k + 1$ & $3k + 2$ & $k + 1$ & 1 \\ \hline
\end{tabular}
\end{center}

\subsection{$c_5 = 0$}.
We have the signatures $\Sigma(-1, 0) = ****\pp+\pp$ and $\Sigma(1, 0) = ****\pp-\pp$
which imply:
\begin{equation*}
\begin{split}
-1 & \leq c_1 \leq 1, \\
-2 & \leq c_2 \leq 2, \\
-3 & \leq c_3 \leq 3, \\
-1 & \leq c_4 \leq 1.
\end{split}
\end{equation*}
As $c_5 = 0$, we can assume without loss of generality $c_1 \geq 0$. We refine by case
distinction by the values of $c_1$.

$c_1 = 0$. Here we can assume without loss of generality that $c_4 \geq 0$.
Let $c_4 = 0$ and thus without loss of generality $c_2 \geq 0$. We have the following
table

\begin{center}
\begin{tabular}{|c|c|c|c|} \hline
$c_2$ & $c_3$ & $m$ & $\Sigma(m)$\\ \hline
0 & 2 & $(-1, -1)$ & $\pp-\pp---+$\\ \hline
0 & 3 & $(-1, -1)$ & $\pp-+---+$\\ \hline
1 & $\leq 0$ & $(1, 0)$ & $-\pp-\pp++-$\\ \hline
1 & $\geq 2$ & $(0, 1)$ & $-\pp+\pp++-$\\ \hline
\end{tabular}
\end{center}
Thus we have found:
\begin{center}
\begin{tabular}{|c|c|c|c|c|} \hline
$c_1$ & $c_2$ & $c_3$ & $c_4$ & $c_5$\\ \hline
0 & 0 & 0 & 0 & 0 \\ \hline
0 & 0 & 1 & 0 & 0 \\ \hline
0 & 1 & 1 & 0 & 0 \\ \hline
\end{tabular}
\end{center}
Now let $c_4 = 1$. Then $\Sigma(-1, -1) = \pp**\pp--+$ and hence $c_2 = 1$ and $c_3 =
2$. We have found:
\begin{center}
\begin{tabular}{|c|c|c|c|c|} \hline
$c_1$ & $c_2$ & $c_3$ & $c_4$ & $c_5$\\ \hline
0 & 1 & 2 & 1 & 0 \\ \hline
\end{tabular}
\end{center}

$c_1 = 1$. Assume first that $c_4 = -1$. Then $\Sigma(0, 1) = \pp**-++-$ which makes
$c_4 = -1$ impossible.

Now let $c_4 = 0$. We have $\Sigma(0, 1) = \pp**\pp++-$ which implies that $c_2 = 1$
and $c_3 = 1$. We have found
\begin{center}
\begin{tabular}{|c|c|c|c|c|} \hline
$c_1$ & $c_2$ & $c_3$ & $c_4$ & $c_5$\\ \hline
1 & 1 & 1 & 0 & 0 \\ \hline
\end{tabular}
\end{center}

Finally, let $c_4 = 1$. Then $\Sigma(0, 0) = +**+\pp\pp\pp$, so
\begin{equation*}
\begin{split}
c_2 & \geq 1, \\
c_3 & \geq 1.
\end{split}
\end{equation*}
So we have reduced to six possibilities. Consider the table
\begin{center}
\begin{tabular}{|c|c|c|c|} \hline
$c_2$ & $c_3$ & $m$ & $\Sigma(m)$\\ \hline
1 & 1 & $(-1, -1)$ & $+\pp-\pp--+$\\ \hline
1 & 3 & $(-1, -1)$ & $+\pp+\pp--+$\\ \hline
2 & 1 & $(0, 2)$ & $-\pp-+++-$\\ \hline
\end{tabular}
\end{center}
The remaining cases are:
\begin{center}
\begin{tabular}{|c|c|c|c|c|} \hline
$c_1$ & $c_2$ & $c_3$ & $c_4$ & $c_5$\\ \hline
1 & 1 & 2 & 1 & 0 \\ \hline
1 & 2 & 2 & 1 & 0 \\ \hline
1 & 2 & 3 & 1 & 0 \\ \hline
\end{tabular}
\end{center}
which finishes the classification.

\section{Table of cohomology-free line bundles and theorem}
\label{listsection}

We represent the classification obtained in the previous section in a table at the end
of this section.
We distinguish three types of line bundles, named by the letters $A$ to $C$, where the
$B$-type bundles form infinite series. For a given cohomology-free bundle $\mathcal{L}$
the table shows the tuple $(c_1, c_2, c_3, c_4, c_5)$ and a list all cohomology-free
bundles $\mathcal{L}'$ which have the property that $H^k(X, \mathcal{L}^* \otimes
\mathcal{L}) = H^k(X, \mathcal{L} \otimes (\mathcal{L}')^*) = 0$ for all $k > 0$, which
is a necessary condition for $\mathcal{L}$ and $\mathcal{L}'$ for being part of the
same strongly exceptional sequence. We say that $\mathcal{L}$ and $\mathcal{L}'$ are
{\em compatible}. For notation, $-A_4$ for instance means the line
bundle $(-1, -1, -1, 0, 0)$. Now we state and proof our main result. Let $X$ be the
toric surface as given in the introduction.

\begin{theorem}
On $X$ there are no strongly exceptional sequences of length 7 which consist of
line bundles.
\end{theorem}

\begin{proof}
The proof is done by inspection of the table and exclusion principle.
For example, assume that
we have a strongly exceptional sequence of length 7 which contains $C_{10}$. Then the
rest of the sequence can at most be selected from $A_1$, $A_2$, $A_4$, $C_3$, $C_7$,
$C_9$, $B_{4, 1}$, $B_{4, 2}$. We see from the corresponding rows that at most one of
the $A_i$ and at most one of the $C_i$ can be selected simultaneously. Hence we can
chose at most four elements from the list to complete the sequence. We conclude that a
strongly exceptional sequence of length 7 which contains $C_{10}$ cannot exist. Thus
we can eliminate $C_{10}$ from the table.

As general rules we read off that at most two of the $A_i$ can be part of a strongly
exceptional sequence, i.e. we have either $\pm A_i$, $i = 1, \dots, 7$ alone or $A_i$,
$i = 1, \dots, 7$, and $A_7$ (respectively $-A_i$ and $-A_7$), together, or $A_i$ and
$-A_{7 - i}$, $i = 1, \dots, 6$ (respectively $-A_i$ and $A_{7 - i}$), together.

Assume that a strongly exceptional sequence contains three bundles of type $B_{r, k}$,
$B_{s, l}$, $B_{t, m}$. We read immediately off from the table that this is not
possible
if $r, s, t$ are pairwise distinct, hence at least two of the $r, s, t$ coincide. We
also see that always $B_{r, k + 1} - B_{r, k} = A_7$ for all $r$ and $B_{r, k + n} -
B_{r, k} = n \cdot A_7$, so if two bundles of the same $B$-type are contained in
a strongly exceptional sequence, these must be of the form $B_{r, k}$, $B_{r, k + 1}$.
Now given such a pair and assume that there exists one more $B_{s, l}$ together with
this pair in a strongly exceptional sequence. Then $B_{r, k + 1} - B_{s, l} = A_i$ for
some $1 \leq i \leq 6$ and $B_{r, k} - B_{s, l} = -A_{7 - i}$. If there exists another
$B_{t, m}$ in this sequence, we have $B_{r, k + 1} - B_{t, m} = A_j$ for $1 \leq i \leq
6$ and $B_{t, m} - B_{s, l} = A_i - A_j$, which is not possible. So we conclude that
a strongly exceptional sequence can contain at most three of the $B$'s.
This in turn, together with the above condition on the $A$'s, implies that a strongly
exceptional sequence must contain at least one of the $C$'s.

We proceed now with $C_9$. We can only choose at most three of the compatible $B$'s
and at most one of the $A$'s. So we have to choose at least one out of $C_2$ and
$C_6$. These two are mutally exclusive, so we can choose only one of them. Both
choices restrict the choice of the $A$'s to $-A_1$. $-A_1$ in turn is not compatible
with $B_{4, k}$, so that we can choose at most two of the $B$'s, which is not enough,
hence we can forget about $C_9$.

For $C_8$, we can choose at mosts three of the $B$'s and thus to obtain a strongly
exceptional sequence, we have to choose both, $A_5$ and $C_1$. But $A_5$ is not
compatible
with $B_2$, so we can not complete to a full sequence. Hence we eliminate $C_8$.

$C_7$. The bundles $C_1$ and $C_6$ are mutually exclusive, so in order to obtain an
exceptional sequence of length seven, we have to choose one
out of the $A$'s and three out of the $B$'s. The $C$'s leave only one choice for the
$A$'s, namely $-A_2$, which in turn is not compatible with $B_{4, k}$, hence we can
discard $C_7$.

$C_6$. Here we have only the choice of at most one of the $A$'s and of at most three
of the $B$'s left, which is not enough. So $C_6$ goes away.

$C_5$. Both pairs $C_3$, $C_4$ and $B_{1, 2}$, $B_{7, 2}$ are mutually exclusive,
leaving not enough choices to complete the sequences. Bye bye, $C_5$.

$C_4$. The sequence must contain $C_1$ and $-A_2$, where the latter is not compatible
with the $B_7$'s, so no $C_4$.

$C_3$. We can choose at most one $A$ and at most one $C$. The $C$'s are not compatible
with $A_1$ and $A_2$, and $B_4$ and $B_7$ are not simultaneously compatible with one
of $-A_3$ and $-A_4$, which does not leave enough choices. to choose also $-A_4$,
which is not compatible with $B_{4, k}$. So we can also exclude $C_3$.

In the remaining cases, for $C_1$ and $C_2$, we do not have any other $C$'s at our
disposal. Therefore we can not complete to a sequence and so we can eliminate
$C_1$ and $C_2$.

Altogether, we have removed now all $C$'s, and as we have seen above, it is not
possible to complete to a strongly exceptional sequence of length 7.
\end{proof}

\begin{center}
\begin{tabular}{|c|c|c|} \hline
Name & $(c_1, c_2, c_3, c_4, c_5)$ & Compatible with \\ \hline \hline
$A_1$ & $(0, 0, 1, 0, 0)$ & $-A_6$, $A_7$, $-C_2$, $C_3$, $-C_6$, $C_7$, $-C_9$,
$C_{10}$\\ 
& & $-B_{1, k}$, $B_{2, k}$, $-B_{3, k}$, $B_{4, k}$, $-B_{6, k}$, $B_{7, k}$ \\ \hline

$A_2$ & $(0, 1, 1, 0, 0)$ & $-A_5$, $A_7$, $-C_1$, $C_3$, $-C_4$, $C_5$, $-C_6$,
$-C_7$, $C_9$, $C_{10}$\\
& & $-B_{1, k}$, $B_{3, k}$, $-B_{2, k}$, $B_{4, k}$, $-B_{5, k}$, $B_{7, k}$ \\ \hline

$A_3$ & $(0, 1, 2, 1, 0)$ & $-A_4$, $A_7$, $-C_1$, $-C_3$, $C_4$, $C_5$\\
& &  $-B_{2, k}$, $B_{5, k}$,  $-B_{3, k}$, $B_{6, k}$,  $-B_{4, k}$, $B_{7, k}$
if $k \geq 1$ \\ \hline

$A_4$ & $(1, 1, 1, 0, 0)$ & $-A_3$, $A_7$, $-C_1$, $-C_2$, $-C_3$, $C_7$, $C_9$,
$C_{10}$, \\
& & $B_{2, k}$ if $k \geq 2$, $B_{3, k}$ if $k \geq 2$, $B_{4, k}$, \\
& & $-B_{5, k}$, $-B_{6, k}$, $B_{4, k}$, $-B_{7, k}$ \\ \hline

$A_5$ & $(1, 1, 2, 1, 0)$ & $-A_2$, $A_7$, $-C_1$, $C_8$, $B_{1, k}$, $B_{2, k}$
if $k \geq 2$, $-B_{3, k}$,\\
& & $-B_{4, k}$, $B_{5, k}$, $-B_{6, k}$ \\ \hline

$A_6$ & $(1, 2, 2, 1, 0)$ & $-A_1$, $A_7$, $B_{1, k}$, $-B_{2, k}$, $B_{3, k}$
if $k \geq 2$, \\
& & $-B_{4, k}$, $B_{6, k}$, $-B_{7, k}$ \\ \hline

$A_7$ & $(1, 2, 3, 1, 0)$ & $A_1$, $A_2$, $A_3$, $A_4$, $A_5$, $A_6$, $B_{1, k}$
for $k \geq 3$, $B_{7, 1}$, \\
& & $B_{i, k}$ for $i = 2, \dots, 7$ and $k \geq 2$, \\
& & $-B_{i, k}$ for
$i = 1, \dots, 7$. \\ \hline \hline

$B_{1, k}$, $k \geq 2$ & $(k, 2k - 1, 3k - 2, k, 1)$ & $-A_1$, $-A_2$, $A_5$,
$A_6$, $A_7$ if $k \geq 3$, $-A_7$, \\
& & $B_{1, k - 1}$ if $k \geq 2$, $B_{1, k + 1}$, $B_{2, k - 1}$, $B_{2, k}$,
$B_{3, k - 1}$, $B_{3, k}$ \\ \hline

$B_{2, k}$, $k \geq 1$ & $(k, 2k - 1, 3k - 1, k, 1)$ & $A_1$, $-A_2$, $-A_3$, $A_4$ if
$k \geq 2$,
$-A_6$, $A_7$ if $k \geq 2$, $-A_7$,\\
& & $B_{1, k}$ if $k \geq 2$, $B_{1, k + 1}$, $B_{2, k - 1}$ if $k \geq 2$, $B_{2, k + 1}$, \\
& & $B_{4, k - 1}$ if $k \geq 2$, $B_{4, k}$,
$B_{5, k - 1}$ if $k \geq 2$, $B_{5, k}$ \\ \hline

$B_{3, k}$, $k \geq 1$ & $(k, 2k, 3k - 1, k, 1)$ & $-A_1$, $A_2$, $A_4$ if $k \geq 2$,
$-A_5$, $A_6$ if $k \geq 2$, $A_7$ if $k \geq 2$, $-A_7$,\\
& & $B_{1, k}$ if $k \geq 2$, $B_{1, k + 1}$, $B_{3, k - 1}$ if $k \geq 2$, $B_{3, k + 1}$ \\
& & $B_{4, k - 1}$ if $k \geq 2$, $B_{4, k}$, $B_{6, k - 1}$ if $k \geq 2$, $B_{6, k}$ \\ \hline

$B_{4, k}$, $k \geq 1$ & $(k, 2k, 3k, k, 1)$ & $A_1$, $A_2$, $-A_3$, $A_4$ if
$k \geq 2$, $-A_5$, $-A_6$, $A_7$ if $k \geq 2$, $-A_7$, \\
& & $B_{2, k}$, $B_{2, k + 1}$,
$B_{3, k}$, $B_{3, k + 1}$, $B_{4, k - 1}$ if $k \geq 2$, $B_{4, k + 1}$, \\
& & $B_{7, k - 1}$ if $k \geq 2$, $B_{7, k}$ \\ \hline

$B_{5, k}$, $k \geq 1$ & $(k, 2k, 3k + 1, k + 1, 1)$ & $-A_2$, $-A_4$, $A_5$ if
$k \geq 2$, $A_7$ if $k \geq 2$, $-A_7$, \\
& & $B_{2, k}$, $B_{2, k + 1}$,
$B_{5, k - 1}$ if $k \geq 2$, $B_{5, k + 1}$,\\
& & $B_{7, k - 1}$ if $k \geq 2$, $B_{7, k}$ \\ \hline

$B_{6, k}$, $k \geq 1$ & $(k, 2k + 1, 3k + 1, k + 1, 1)$ & $-A_1$, $A_3$, $-A_4$,
$-A_5$, $A_6$ if $k \geq 2$, $A_7$ if $k \geq 2$, $-A_7$, \\
& & $B_{3, k}$, $B_{3, k + 1}$,
$B_{6, k - 1}$ if $k \geq 2$, $B_{6, k + 1}$,\\
& & $B_{7, k - 1}$ if $k \geq 2$, $B_{7, k}$ \\ \hline

$B_{7, k}$, $k \geq 0$ & $(k, 2k + 1, 3k + 2, k + 1, 1)$ & $A_1$ if $k \geq 1$,
$A_2$ if $k \geq 1$, $A_3$ if $k \geq 1$,
$-A_4$, $-A_5$,\\
& & $-A_6$, $A_7$ if $k \geq 1$, $-A_7$, \\
& & $B_{4, k}$, $B_{4, k + 1}$,
$B_{5, k}$, $B_{5, k + 1}$, $B_{6, k}$, $B_{6, k + 1}$,\\
& & $B_{7, k - 1}$ if $k \geq 2$, $B_{7, k + 1}$ \\ \hline
\end{tabular}
\end{center}

\begin{center}
\begin{tabular}{|c|c|c|} \hline
Name & $(c_1, c_2, c_3, c_4, c_5)$ & Compatible with \\ \hline \hline
$C_1$ & $(2, 4, 7, 3, 2)$ & $-A_2$, $-A_3$, $-A_4$, $-A_5$, $C_3$, $C_4$, $C_7$,
$C_8$, \\
& & $B_{2, 1}$, $B_{2, 2}$, $B_{4, 1}$, $B_{5, 1}$, $B_{7, 0}$, $B_{7, 1}$\\ \hline

$C_2$ & $(2, 5, 7, 3, 2)$ & $-A_1$, $-A_4$, $C_3$, $C_9$, \\
& & $B_{3, 1}$, $B_{3, 2}$, $B_{4, 1}$, $B_{6, 1}$, $B_{7, 0}$, $B_{7, 1}$ \\ \hline

$C_3$ & $(2, 5, 8, 3, 2)$ & $A_1$, $A_2$, $-A_3$, $-A_4$, $C_1$, $C_2$,
$C_5$, $C_{10}$, \\
& & $B_{4, 1}$, $B_{4, 2}$, $B_{7, 0}$, $B_{7, 1}$ \\ \hline

$C_4$ & $(2, 5, 9, 4, 2)$ & $-A_2$, $A_3$, $C_1$, $C_5$,\\
& & $B_{5, 1}$, $B_{5, 2}$, $B_{7, 0}$, $B_{7, 1}$ \\ \hline

$C_5$ & $(2, 6, 10, 4, 2)$ & $A_2$, $A_3$, $C_3$, $C_4$, $B_{7, 0}$,
$B_{7, 2}$ \\ \hline

$C_6$ & $(3, 5, 7, 3, 2)$ & $-A_1$, $-A_2$, $C_7$, $C_9$, \\
& & $B_{1, 2}$, $B_{2, 1}$, $B_{2, 2}$, $B_{3, 1}$, $B_{3, 2}$, $B_{4, 1}$ \\ \hline

$C_7$ & $(3, 5, 8, 3, 2)$ & $A_1$, $-A_2$, $A_4$, $C_1$, $C_6$, $C_{10}$ \\
& & $B_{2, 1}$, $B_{2, 2}$, $B_{4, 1}$, $B_{4, 2}$\\ \hline

$C_8$ & $(3, 5, 9, 4, 2)$ & $A_5$, $C_1$, \\
& & $B_{2, 1}$, $B_{2, 2}$, $B_{5, 1}$, $B_{5, 2}$ \\ \hline

$C_9$ & $(3, 6, 8, 3, 2)$ & $-A_1$, $A_2$, $A_4$, $C_2$, $C_6$, $C_{10}$ \\
& & $B_{3, 1}$, $B_{3, 2}$, $B_{4, 1}$, $B_{4, 2}$ \\ \hline

$C_{10}$ & $(3, 6, 9, 3, 2)$ & $A_1$, $A_2$, $A_4$, $C_3$, $C_7$, $C_9$ \\
& & $B_{4, 1}$, $B_{4, 2}$ \\ \hline
\end{tabular}
\end{center}

\end{document}